\documentclass[a4in,12pt]{article}
\usepackage{amsmath}
\usepackage{amssymb}
\usepackage{amscd}

\newtheorem{theoreme}{Th\'eor\`eme}[section]
\newtheorem{proposition}[theoreme]{Proposition}
\newtheorem{corollaire}[theoreme]{Corollaire}
\newtheorem{lemme}[theoreme]{Lemme}

\newtheorem{exemple}[theoreme]{Exemple}
\newtheorem{remarque}[theoreme]{Remarque}

\newenvironment{preuve}{\begin{trivlist} \item[]{\it Preuve---}}
{\par\hfill $\square$\end{trivlist}}
\newcommand{\Aut}{{\rm Aut}}
\newcommand{\rank}{{\rm rank}}

\renewcommand{\P}{\mathbb{P}}
\newcommand{\PP}{{\cal P}}
\newcommand{\C}{\mathbb{C}}

\newcommand{\R}{\mathbb{R}}

\newcommand{\Z}{\mathbb{Z}}
\newcommand{\id}{{\rm id}}
\newcommand{\h}{{\rm h}}

\renewcommand{\H}{{\cal H}}
\newcommand{\HH}{{\rm H}}
\newcommand{\G}{{\cal G}}
\newcommand{\K}{{\cal K}}

\newcommand{\cad}{{\it c.-\`a-d. }}
\newcommand{\ie}{{\it i.e. }}

\newcommand{\dist}{{\rm dist}}

\newcommand{\F}{{\cal F}}

\newcommand{\E}{{\cal E}}

\newcommand{\voir}{{\it voir }}

\title{Groupes commutatifs d'automorphismes d'une vari\'et\'e
  k\"ahl\'erienne compacte}
\author{Tien-Cuong Dinh et Nessim Sibony}
\begin{document}
\maketitle
\begin{abstract}
Let $V$ be a compact K\"ahler manifold. Let $\G'$
be a commutative subgroup of $\Aut(V)$ and $U$ the set of elements of 
zero entropy of $\G'$.
Then $U$ is a group and $\G'$ is isomorphic to the direct product
$\G'\simeq U\times \G$ where $\G$ is a subgroup of $\G'$ such that all 
elements of $\G\setminus\{\id\}$ are of positive entropy. 
Moreover, $\G$ is a free commutative group with $\rank(\G)\leq \dim V-1$. 
The estimate is sharp. When $\rank(\G)=\dim V-1$, $U$ is finite. 
$\mbox{Rank}(\G)$ satisfies other inequalities involving the dimensions of
Dolbeault cohomology groups of $V$.
\end{abstract}
\section{Introduction}
Soit $f$ un endomorphisme holomorphe d'une vari\'et\'e $V$
k\"ahl\'erienne compacte. On peut consid\'erer qu'un tel objet est
dynamiquement int\'eressant lorsqu'il est d'entropie positive,
c'est-\`a-dire lorsque le nombre d'orbites distinctes discernables \`a
l'\'echelle $\epsilon$ est \`a croissance exponentielle. D'apr\`es 
Gromov
\cite{Gromov3} et Yomdin \cite{Yomdin} (\voir \'egalement
Katok-Hasselblatt \cite{KatokHasselblatt}),
un endomorphisme 
holomorphe $f$ est d'entropie positive si le rayon spectral de
$f^*$ sur $\H^{1,1}(V,\R)$ est strictement plus grand que $1$, 
$\H^{1,1}(V,\R)$ d\'esignant le groupe de cohomologie de Dolbeault.
\par
On trouve dans la litt\'erature des exemples d'automorphismes
d'entropie positive; pour le cas des surfaces \voir B. Mazur \cite{Mazur}, S.
Cantat \cite{Cantat1}, Barth, Peters, van de Ven \cite{BarthPetersVen}.
La construction de Mazur se prolonge sans difficult\'e en dimension $k$.
Tr\`es r\'ecemment, C. McMullen \cite{McMullen} a construit des
surfaces K3 non alg\'ebriques ayant des automorphismes admettant des disques de
Siegel. Il d\'emontre en particulier que ces automorphismes sont d'entropie
positive. En fait, on doit observer qu'il est assez difficile de construire des 
automorphismes holomorphes de vari\'et\'es complexes. C'est aussi le cas dans le 
cadre de domains strictement pseudoconvexes, \voir Burns-Shnider-Wells
\cite{BurnsShniderWells}.  
\par
Dans \cite{Smale}, S. Smale attire l'attention sur la
difficult\'e de comprendre les diff\'eomorphismes d'une
vari\'et\'e qui commutent.
Dans le pr\'esent article, nous explorons les propri\'et\'es des
groupes commutatifs  d'automorphismes d'entropie positive. Plus
pr\'esis\'ement, un groupe  $\G\subset \Aut(V)$ est d'entropie
positive si tout \'el\'ement $g\in\G\setminus\{\id\}$ est
d'entropie positive. Nous montrons
\vspace{0.3cm}
\\
{\bf Th\'eor\`eme principal.} {\it Si $\G$ est un groupe commutatif
d'entropie positive d'une vari\'et\'e K\"ahl\'erienne
compacte $V$. Alors $\G$ est commutatif 
libre et le rang $r$ de $\G$ v\'erifie 
$r\leq \dim V-1$.}
\vspace{0.3cm}
\\
La preuve utilise de fa\c con essentielle une extension du
th\'eor\`eme de Hodge-Riemann due \`a Gromov \cite{Gromov2} 
(\voir aussi Timorin \cite{Timorin}). Elle r\'esulte de
l'analyse de l'action de $\G$ sur $\H^{1,1}(V,\R)$. Notons que
nous montrons \`a l'aide d'exemples classiques (\voir exemple 4.7),
fournis par F. Labourie et F. Paulin, que notre estimation est
optimale.
Nous en d\'eduisons gr\^ace \`a un th\'eor\`eme de Demailly-Paun sur le c\^one 
de K\"ahler,  la structure des groupes 
commutatifs quelconques. 
\par
L'analyse des endomorphismes holomorphes qui commutent dans
$\P^k$ a fait l'objet d'un pr\'ec\'edent travail des auteurs \cite{DinhSibony1}. 
Le
cas $\P^1$ a \'et\'e \'etudi\'e par Fatou \cite{Fatou}, Julia
\cite{Julia}, Ritt \cite{Ritt} et Eremenko \cite{Eremenko}. Dans
le cadre non compact signalons que Veselov \cite{Veselov}
et Lamy \cite{Lamy} ont montr\'e que tout groupe d'automorphismes
polynomiaux de $\C^2$ satisfait \`a l'alternative de Tits,
c'est-\`a-dire que soit il est moyennable soit il contient un sous
groupe libre \`a 2 g\'en\'erateurs. 
\par
Nous remercions E. Amerik, Y. Benoist, J. B. Bost, 
F. Labourie et F. Paulin d'avoir 
r\'epondu \`a nos nombreuses questions.
\section{Action sur les groupes de cohomologie}
Soit $V$ une vari\'et\'e k\"ahl\'erienne compacte de dimension $k$
munie d'une forme de K\"ahler $\omega$.
Un endomorphisme (non d\'eg\'en\'er\'e) de $V$ 
induit des applications lin\'eaires inversibles sur les groupes de
cohomologie de la vari\'et\'e. 
Nous allons rappeler la d\'efinition 
des degr\'es dynamiques et de l'entropie associ\'es \`a un
endomorphisme $f$ de $V$.
\par
Notons $\F_\C^{p,p}$ (resp. $\E_\C^{p,p}$) l'espace des
$(p,p)$-formes lisses, $\overline\partial$-ferm\'ees
(resp. $\overline\partial$-exactes) de
$V$ et $\F_\R^{p,p}$ (resp. $\E_\R^{p,p}$) le sous-espace des formes 
$\psi\in \F_\C^{p,p}$ (resp. $\psi\in \E_\C^{p,p}$)
v\'erifiant $\psi=\overline\psi$. 
Consid\'erons les groupes de cohomologie de Dolbeault \cite{GriffithsHarris}
$$\H^{p,p}(V,\C):=\F_\C^{p,p}/\E_\C^{p,p} \mbox{ \ \ et \ \ }
\H^{p,p}(V,\R):=\F_\R^{p,p}/\E_\R^{p,p}.$$
Ces espaces sont de dimension
finie et on a $\H^{p,p}(V,\C)=\H^{p,p}(V,\R)\otimes_\R\C$. On peut
identifier  $\H^{p,p}(V,\R)$ \`a un sous-espace r\'eel de
$\H^{p,p}(V,\C)$. 
\par 
Si dans la d\'efinition pr\'ec\'edente, on remplace les
formes par les courants, on obtient les m\^emes groupes
de cohomologie. C'est une cons\'equence de la dualit\'e de
Poincar\'e et de la d\'ecomposition de Hodge.
Etant donn\'e une
forme ou un courant $\alpha$ ferm\'e 
de bidegr\'e $(p,p)$, on notera $[\alpha]$
sa classe dans $\H^{p,p}(V,\C)$.
Notons $\K$ {\it le c\^one de K\"ahler} de 
$\H^{1,1}(V,\R)$. C'est le c\^one des classes associ\'ees aux formes
de K\"ahler, \ie $(1,1)$-formes lisses strictement positives.
C'est un c\^one ouvert de $\H^{1,1}(V,\R)$. 
Observons que 
$\overline \K$ est convexe et saillant: $\overline \K\cap (-\overline \K)=\{0\}$.
En effet, chaque classe de $\overline \K$ 
est repr\'esent\'ee par un courant positif
ferm\'e.
\par
Consid\'erons un endomorphisme holomorphe $f$ de $V$.
L'application $f^*$ d\'efinie par $f^*[\alpha]:=[f^*\alpha]$, 
est un op\'erateur
lin\'eaire de $\H^{p,p}(V,\R)$ dans lui-m\^eme. 
Posons
$$d_{p,n}:= 
\|(f^n)^*[\omega^p]\|$$
$$d_p:=\limsup_{n\rightarrow\infty}
\sqrt[n]{d_{p,n}} \mbox{ \ et \ }
\HH(f):=\sup_{0\leq p\leq k}\log (d_p)$$
o\`u $\|\ \|$ d\'esigne une norme sur $\H^{p,p}(V,\R)$. Observons
que $d_{p,n}$ ne d\'epend que de la classe $[\omega^p]$
dans $\H^{p,p}(V,\R)$. Il est aussi
clair que $d_p$ et $\HH(f)$ ne d\'ependent pas de la forme
$\omega$ et de la norme $\|\ \|$ choisies. De plus, $d_p$ est le rayon
spectral de $f^*$ sur $\H^{p,p}(V,\R)$ et sur  $\H^{p,p}(V,\C)$.
On dit que $d_p$ est le degr\'e dynamique d'ordre $p$ de $f$.
En pratique, on peut calculer $d_p$ par la formule 
$$d_p:=\lim_{n\rightarrow \infty} \left(\int (f^n)^*\omega^p\wedge
\omega^{k-p} \right)^{1/n}.$$
\par
Rappelons la notion d'entropie topologique \cite{KatokHasselblatt}. 
Notons $\dist$
la distance sur $V$. On d\'efinit la distance $\dist_n$ par 
$$\dist_n(a,b):=\max_{0\leq i\leq n-1} \dist(f^i(a), f^i(b)).$$
Soit $s_n(\epsilon)$ le nombre maximal de boules de rayon $\epsilon/2$
pour la distance $\dist_n$ qui sont disjointes dans $V$. {\it L'entropie
  topologique} de $f$ est d\'efinie par la formule 
$$\h(f):= \limsup_{\epsilon\rightarrow
  0}\limsup_{n\rightarrow\infty}\frac{\log s_n(\epsilon)}{n}.$$
Nous avons le r\'esultat fondamental suivant
\cite{Yomdin,Gromov3, Friedland}. 
\begin{theoreme}[Yomdin-Gromov] Soit $f$ un endomorphisme
holomorphe d'une
vari\'et\'e k\"ahl\'erienne compacte $V$. 
Alors l'entropie topologique
$\h(f)$ de $f$
est \'egale \`a $\HH(f)$. 
\end{theoreme} 
\par
On en d\'eduit le corollaire suivant.
\begin{corollaire} On a
 $\h(f)>0$ si et seulement si $d_1>1$. De plus,
l'application qui associe \`a un endomorphisme holomorphe de $V$
  son premier degr\'e dynamique $d_1$ est d'image discr\`ete dans
  $[1,+\infty[$.  
L'application qui associe \`a un endomorphisme holomorphe de $V$
  son entropie topologique est d'image discr\`ete dans
  $[0,+\infty[$.  
\end{corollaire}
\begin{preuve} Soit $f$ un endomorphisme holomorphe de $V$.
L'int\'egrale d\'efinissant $d_{p,n}$ se calcule cohomologiquement. Il
  suffit d'utiliser la d\'ecomposition de Jordan de $f^*$ sur
  $\H^{1,1}(V,\C)$ pour voir que  
  $d_p$ est \'egal au module d'un produit de $p$ valeurs propres de
  $f^*$ sur $\H^{1,1}(V,\C)$. En effet, la norme se calcule sur
  $\omega^p$ qui est un produit de $(1,1)$ formes.
En particulier, on a $d_p\leq d_1^p$. D'apr\`es le th\'eor\`eme 2.1,
il nous suffit de montrer que lorsque 
  le premier degr\'e dynamique $d_1$ de $f$ est 
major\'e par une constante $M>0$, les valeurs
  propres de $f^*$ sur le groupe de cohomologie de De Rham
$\H^2(V,\C)$ parcourent un ensemble fini. 
\par
La d\'ecomposition de Hodge 
$$\H^2(V,\C)=\H^{2,0}(V,\C)\oplus \H^{1,1}(V,\C)\oplus \H^{0,2}(V,\C)$$
est invariante par $f^*$.
Si $\psi$ est une $(2,0)$-forme holomorphe non nulle telle que
$f^*\psi=\lambda\psi$ avec $\lambda\in \C$, on a
$$f^*(\psi\wedge\overline \psi)=|\lambda|^2 \psi\wedge\overline
\psi.$$
La forme $\psi\wedge\overline\psi$ \'etant lisse, positive et non nulle, il
existe $c>0$ telle que $\psi\wedge\overline\psi\leq c\omega$. 
On en d\'eduit que
\begin{eqnarray*}
|\lambda|^2 & = & \lim_{n\rightarrow\infty}\left( \int
(f^n)^*(\psi\wedge\overline \psi)\wedge \omega^{k-2}\right)^{1/n}\\
& \leq & 
 \lim_{n\rightarrow\infty}\left( \int
c(f^n)^*\omega^2\wedge \omega^{k-2}\right)^{1/n}= d_2.
\end{eqnarray*}
Par cons\'equent, on a 
$|\lambda|^2\leq d_1^2\leq M^2$. On en d\'eduit que le rayon
spectral de $f^*$ sur $\H^{2,0}(V,\C)$ est major\'e par $M$. Par
conjugaison, le rayon
spectral de $f^*$ sur $\H^{0,2}(V,\C)$ l'est aussi. 
\par
Par suite, le rayon spectral de $f^*$ sur $\H^2(V,\C)$ est major\'e par
$M$. D'autre part, l'op\'erateur $f^*$ pr\'eserve $\H^2(V,\Z)$ et  
on a $\H^2(V,\C)=\H^2(V,\Z)\otimes_\Z \C$. 
L'op\'erateur $f^*$ pr\'eserve un r\'eseau qui engendre $\H^2(V,\C)$. 
Par cons\'equent, pour
une base convenable,
il est d\'efini par une matrice \`a coefficients
entiers. Son polyn\^ome caract\'eristique est unitaire et \`a
coefficients entiers. Les racines de ce polyn\^ome sont major\'ees par
$M$. Il n'y a donc qu'un nombre fini de polyn\^omes v\'erifiant
ces propri\'et\'es et par cons\'equent un nombre fini de valeurs
propres possibles pour $f^*$.  
\end{preuve}
\begin{remarque} \rm
De la m\^eme mani\`ere, on montre qu'il existe un nombre entier naturel 
$m\geq 1$ tel que si $f$ est un automorphisme d'entropie nulle de $V$, 
les valeurs propres de $f^{m*}$ sur $\H^{1,1}(V,\R)$ sont \'egales \`a 1.
\end{remarque}
\section{Th\'eor\`eme de Hodge-Riemann}
Nous rappelons dans ce paragraphe le th\'eor\`eme de signature de
Hodge-Riemann. Nous en donnons des corollaires qui nous seront utiles.
\par
Soit $V$ une vari\'et\'e k\"ahl\'erienne compacte de
  dimension $k$. Soient
  $c_1$, $\ldots$, $c_{k-1}$ des classes de  $\H^{1,1}(V,\R)$. 
On d\'efinit la
  forme sym\'etrique $q$ sur $\H^{1,1}(V,\R)$  par
$$q(c,c'):=-\int_V c\wedge c'\wedge c_1\wedge\ldots\wedge c_{k-2}.$$
Notons
$\PP(c_1\wedge \ldots \wedge c_{k-1})$ 
l'ensemble des {\it classes primitives} 
dans $\H^{1,1}(V,\R)$. Plus pr\'ecis\'ement
$$\PP(c_1\wedge \ldots \wedge c_{k-1}):=\big\{c\in \H^{1,1}(V,\R),\ c\wedge
c_1\wedge\ldots\wedge c_{k-1}=0\big\}.$$
Lorsque $c_1\wedge\ldots\wedge c_{k-1}\not =0$, par
dualit\'e de Poincar\'e, $\PP(c_1\wedge\ldots\wedge c_{k-1})$ est un hyperplan
de $\H^{1,1}(V,\R)$.
\par
Rappelons le th\'eor\`eme de Hodge-Riemann classique \cite{GriffithsHarris}:
\begin{theoreme}[Hodge-Riemann] Soit $\omega$ une forme de K\"ahler. Si
  $c_1=\cdots=c_{k-1}=[\omega]$, la forme $q$ est d\'efinie positive
  sur l'hyperplan $\PP([\omega^{k-1}])$.
\end{theoreme}
\begin{corollaire} Soient $c$ et $c'$ deux classes de
  $\overline\K$. Supposons que $c\wedge c'=0$. Alors $c$ et $c'$ sont colin\'eaires.
\end{corollaire}
\begin{preuve} Soit $\omega$ une forme de K\"ahler.
Notons $\PP:=\PP([\omega^{k-1}])$. La forme $q$ est
  d\'efinie positive sur $\PP$. Notons $F$ le sous-espace engendr\'e par
  $c$ et $c'$. Par hypoth\`ese, 
la forme $q$ est d\'efinie semi-n\'egative sur $F$. Par cons\'equent,
  l'intersection de $F$ avec l'hyperplan 
$\PP$ est r\'eduit \`a $\{0\}$. Par suite,
  $\dim F\leq 1$ et donc $c$, $c'$ sont colin\'eaires. 
\end{preuve}
\par
La version suivante du th\'eor\`eme de Hodge-Riemann est d\^ue \`a
Gromov \cite{Gromov2} (\voir aussi Timorin \cite{Timorin}):
\begin{theoreme} Si $c_1$, $\ldots$, $c_{k-1}$ sont des classes de
  K\"ahler, la forme $q$ est d\'efinie semi-positive sur
  $\PP(c_1\wedge \ldots\wedge c_{k-1})$. 
\end{theoreme} 
\par
Par passage \`a limite, on obtient le corollaire suivant:
\begin{corollaire} Supposons que les classes $c_1$, $\ldots$, $c_{k-1}$
  soient dans $\overline \K$. Si
  $c_1\wedge\ldots \wedge c_{k-1}$ est non nulle, 
alors la forme $q$ est 
  semi-positive sur $\PP(c_1\wedge \ldots \wedge c_{k-1})$.
\end{corollaire}
\par
Nous aurons besoin de l'extension suivante du corollaire 3.2.
\begin{corollaire} Soient $c_1$, $\ldots$, $c_m$, $c$ et $c'$ des classes 
dans $\overline \K$, $0\leq m\leq k-2$. 
Supposons que $c\wedge c'\wedge c_1\wedge \ldots\wedge c_m=0$. Alors
  il existe $a$, $b$ r\'eels, $(a,b)\not=(0,0)$, tels que
  $$(ac+bc')\wedge c_1\wedge\ldots\wedge c_m\wedge
  c_1'\wedge\ldots \wedge c_{k-m-2}'=0$$ 
pour toutes les classes 
$c_1'$, $\ldots$, $c_{k-m-2}'$ dans $\H^{1,1}(V,\R)$. Si
$c\wedge c_1\wedge \ldots \wedge c_m$ est non nulle, le couple $(a,b)$ est
  unique \`a une constante multiplicative pr\`es.
\end{corollaire}
\par
On peut supposer que
$c$ et $c'$ ne sont pas colin\'eaires et que
$c\wedge c_1\wedge\ldots\wedge c_m \not =0$, 
$c'\wedge c_1\wedge\ldots\wedge c_m \not =0$. En effet, dans le cas
contraire, le corollaire est \'evident.
Notons $F$ le plan engendr\'e par $c$ et
$c'$. Nous aurons besoin de quelques r\'esultats pr\'eliminaires. 
\begin{lemme}
Soient $c_1'$, $\ldots$, $c_{k-m-2}'$ des classes de K\"ahler. Il
existe $\tilde c \in F$, $\tilde c\not =0$, unique \`a une constante
multiplicative pr\`es, telle que 
$$\tilde c\wedge c_1\wedge \ldots \wedge c_m\wedge c_1'\wedge \ldots
\wedge c_{k-m-2}'=0.$$
\end{lemme}
\begin{preuve}
Fixons une classe de K\"ahler $c_{k-m-1}'$. 
Les classes $c_1$, $\ldots$, $c_m$ \'etant dans l'adh\'erence du c\^one de
K\"ahler et $c_1\wedge\ldots \wedge c_m$ \'etant non nulle, on a 
$$c_1\wedge\ldots\wedge c_m\wedge
  c_1'\wedge\ldots \wedge c_{k-m-1}' \not =0.$$ 
En effet, il existe un courant positif ferm\'e non nul qui
  repr\'esente la classe $c_1\wedge\ldots\wedge c_m$. 
Par cons\'equent, 
$$\PP:=\PP(c_1\wedge\ldots\wedge c_m\wedge
  c_1'\wedge\ldots \wedge c_{k-m-1}')$$
est un hyperplan de $\H^{1,1}(V,\R)$.
\par
Notons $q$ la forme
  sym\'etrique d\'efinie par les $k-2$ classes $c_1$, $\ldots$, $c_m$, $c_1'$,
  $\ldots$, $c_{k-m-2}'$.
D'apr\`es le corollaire
  3.4, $q$ est semi-positive sur $\PP$.
Les classes $c$, $c'$ \'etant dans l'adh\'erence du
  c\^one de K\"ahler et $c\wedge c'\wedge c_1\wedge \ldots\wedge
  c_{k-m-2}'$ \'etant nulle, la forme $q$ est 
  semi-n\'egative sur $F$.
On en d\'eduit que $q$ est nulle sur $F\cap \PP$. L'intersection 
$F\cap \PP$ n'est pas r\'eduit \`a $\{0\}$ car $\PP$ est un hyperplan. 
Il existe donc une classe non nulle $\tilde c\in F\cap \PP$.
On a $q(\tilde c, \tilde c)=0$.
\par 
Comme $q$ est d\'efinie semi-positive sur $\PP$, 
d'apr\`es l'in\'egalit\'e de Cauchy-Schwarz, pour toute $\gamma\in
\PP$, on a
$$|q(\tilde c,\gamma)|^2\leq q(\tilde c,\tilde c)
 q(\gamma,\gamma)=0.$$
Donc $q(\tilde c,\gamma)=0$.
D'autre part, $q(\tilde c,c'_{k-m-1})=0$ car
  $\tilde c$ est primitive. La classe $c'_{k-m-1}$ \'etant non
  primitive, $q(\tilde c,\gamma)=0$ pour toute $\gamma\in
  \H^{1,1}(V,\R)$. Par dualit\'e de Poincar\'e, on a
$$\tilde c\wedge c_1\wedge \ldots\wedge c_m\wedge c_1'\wedge \ldots
\wedge c_{k-m-2}' =0.$$
\par
Montrons l'unicit\'e de $\tilde c$. Si $\tilde c$ n'\'etait pas
unique \`a une constante multiplicative pr\`es, la derni\`ere
in\'egalit\'e serait vraie pour toute classe $\tilde c\in F$. En particulier,
elle est vraie pour $\tilde c= c$. C'est la contradiction recherch\'ee
car on a suppos\'ee $c\wedge c_1\wedge \ldots \wedge c_m\not =0$. 
\end{preuve}
\par
Posons $E:=\H^{1,1}(V,\R)$. Notons $\Sigma$ l'ensemble des $(\tilde
c,c_1',\ldots, c_{k-m-2}')\in F\times E^{k-m-2}$ qui v\'erifient 
$$\tilde c\wedge c_1\wedge \ldots\wedge c_m\wedge c_1'\wedge \ldots
\wedge c_{k-m-2}' =0.$$
C'est un sous-ensemble alg\'ebrique r\'eel de $F\times E^{k-m-2}$. Il est
d\'efini par une famille $\F$ de polyn\^omes homog\`enes de degr\'e $k-m-1$
qui sont lin\'eaires en chaque variable $\tilde c$, $c_1'$, $\ldots$,
$c_{k-m-2}'$. D'apr\`es le lemme 3.6, $\Sigma\cap F\times
\K^{k-m-2}$ est une hypersurface irr\'eductible. C'est en fait un graphe au dessus de $\K$. 
Par cons\'equent, les polyn\^omes de $\F$ admettent un facteur commun
non constant. Soit $Q$ le facteur commun irr\'eductible qui s'annulle
sur $\Sigma\cap F\times \K^{k-m-2}$. Notons $\Sigma^*$ l'ensemble des
z\'eros de $Q$. On a $\Sigma^*\subset \Sigma$.
Observons que $\Sigma^*$ est invariant par toute
permutation de $c_1'$, $\ldots$, $c_{k-m-2}'$. En effet, les \'equations d\'efinissant $\Sigma$ sont clairement sym\'etriques et au dessus de $\K^{k-m-2}$, on a vu que $\Sigma$ est un graphe donc $\Sigma^*$ est aussi un graphe.
\begin{lemme} Pour toutes classes de K\"ahler $c_1'$, $\ldots$,
  $c_{k-m-3}'$, il existe une classe non nulle $\tilde c\in F$ telle
  que 
 $$\tilde c\wedge c_1\wedge \ldots\wedge c_m\wedge c_1'\wedge \ldots
\wedge c_{k-m-2}' =0$$
pour toute classe $c_{k-m-2}'\in E$. 
\end{lemme} 
\begin{preuve} 
Posons $\tilde Q(\tilde c, c_{k-m-2}') := Q(\tilde c,
  c_1',\ldots, c_{k-m-2}')$ et 
$$\Sigma_{\tilde c} :=\Big\{c_{k-m-2}'\in E, \ \tilde Q(\tilde c,
c_{k-m-2}')=0\Big\}.$$
Supposons qu'aucune classe non nulle $\tilde c\in F$ ne v\'erifie le lemme. 
On a $\tilde Q(\tilde c,.)\not =0$ pour toute classe $\tilde c\not =0$. Par
cons\'equent, d'apr\`es le lemme 3.6, pour tout $\tilde c\not=0$, 
l'ensemble $\Sigma_{\tilde c}$ est un hyperplan de $E$. 
\par
Notons $h_{\tilde c}$ la forme sym\'etrique d\'efinie par 
$$h_{\tilde c} (\alpha,\gamma):=\int \tilde c\wedge c_1\wedge \ldots
\wedge c_m \wedge c_1'\wedge \ldots \wedge c_{k-m-3}' \wedge
\alpha\wedge \gamma$$
pour $\alpha$ et $\gamma$ dans $E$. On a $h_{\tilde
  c}(\alpha,\gamma)=0$ pour toutes classes 
$\alpha\in \Sigma_{\tilde c}$ et
$\gamma\in E$. 
\par
La forme $h_{\tilde c}$ est nulle sur l'hyperplan 
$\Sigma_{\tilde c}$ de $E$ sans \^etre identiquement nulle. 
Elle est donc d\'efinie semi-positive ou
semi-n\'egative. La signature de $h_{\tilde c}$ 
d\'ependant contin\^ument de $\tilde c$, elle est donc constante. 
Or les signatures de $h_{\tilde c}$ et $h_{-\tilde c}$ sont
oppos\'ees. C'est la contradiction cherch\'ee.
\end{preuve}
{\it D\'emonstration du corollaire 3.5 ---} 
Notons $\tau: F\times E^{k-m-2} \longrightarrow F\times E^{k-m-3}$ la
projection d\'efinie par
$$\tau(\tilde c, c_1',\ldots, c_{k-m-2}') := (\tilde c, c_1', \ldots ,
  c_{k-m-3}').$$  
Posons 
$$\tilde \Sigma := \Big\{(\tilde c, c_1', \ldots, c_{k-m-3}')\in F\times
  E^{k-m-3},\ \tau^{-1}(\tilde c, c_1', \ldots, c_{k-m-3}')\subset
  \Sigma^* \Big\}.$$
C'est un sous-ensemble alg\'ebrique r\'eel de $F\times E^{k-m-3}$. D'apr\`es
  le lemme 3.7, $\tilde \Sigma$ est de codimension 1. 
On en d\'eduit que le polyn\^ome $Q$ contient un facteur qui ne
  d\'epend pas de $c_{k-m-2}'$. Comme $Q$ est
  est irr\'eductible, il est 
  ind\'ependant de $c_{k-m-2}'$. L'ensemble des z\'eros de $Q$ \'etant
  invariant par les
  permutations de $c_1'$, $\ldots$, $c_{k-m-2}'$, ce polyn\^ome $Q$ est
  \'egalement ind\'ependant de $c_1'$, $\ldots$, $c_{k-m-3}'$. Ceci
  implique le corollaire. La classe $\tilde c$ est unique \`a une
  constante multiplicative pr\`es car sinon $Q$ serait identiquement
  nul. Les assertions sur les ensembles alg\'ebriques r\'eels se
  v\'erifient en complexifiant.
\par
\hfill $\square$
\begin{remarque} \rm
Si dans le corollaire 3.5, 
on suppose de plus que $c\wedge c\wedge c_1\wedge \ldots \wedge c_m=0$ et 
$c'\wedge c'\wedge c_1\wedge \ldots \wedge c_m=0$, on n'a plus 
besoin de l'hypoth\`ese 
que les classes $c$ et $c'$ sont dans $\overline \K$.
\end{remarque}
\section{Groupes d'automorphismes}
Dans ce paragraphe, nous \'etudions les groupes commutatifs 
d'automorphismes d'une
vari\'et\'e k\"ahl\'erienne compacte $V$.
Bochner et Montgomery \cite{BochnerMontgomery} ont montr\'e 
que le groupe d'automorphismes holomorphes $\mbox{Aut}(V)$ 
de $V$ 
est un groupe de Lie complexe
de dimension finie pour toute vari\'et\'e complexe compacte $V$. 
\par
Notons $\mbox{Aut}_0(V)$
la composante de l'identit\'e dans $\mbox{Aut}(V)$.
Il est clair que $\mbox{Aut}_0(V)$ agit
trivialement sur les groupes $\H^{p,p}(V,\R)$ et est un sous-groupe 
distingu\'e de  $\mbox{Aut}(V)$; 
le quotient
$\mbox{Aut}^\#(V)$ 
de  $\mbox{Aut}(V)$ par
$\mbox{Aut}_0(V)$ est donc un groupe. Lieberman \cite[Theorem 3.3]{Lieberman},
a montr\'e que $\Aut_0(V)$ est compactifiable. C'est-\`a-dire qu'il existe 
une vari\'et\'e compacte $X$ et une application m\'eromorphe 
$\Phi:X\times V\longrightarrow V$ telles que
$\Aut_0(V)$ peut-\^etre identifi\'e 
\`a un ouvert de Zariski $Y$ de $X$ et l'orbite de tout point $x\in V$ par
$\Aut_0(V)$ est  
\'egal \`a  $\Phi(Y\times\{x\})$.
\par
On dira qu'un automorphisme $f$ est {\it virtuellement
isotope \`a l'identit\'e} 
si $f$ poss\`ede des it\'er\'es appartenant
\`a $\mbox{Aut}_0(V)$. Si $c$ est une classe dans le c\^one de K\"ahler, 
notons $\Aut_c(V)$ le groupe des automorphismes $f$ qui pr\'eservent $c$: 
$f^*c=c$.
D'apr\`es Lieberman \cite[Proposition 2.2]{Lieberman}, 
$\Aut_c(V)/\Aut_0(V)$ est fini. En particulier, les \'el\'ements de 
$\Aut_c(V)$ sont virtuellement isotopes \`a l'identit\'e. Si $V$ n'admet pas
de champ de vecteurs holomorphe, le groupe $\Aut_c(V)$ est fini. 
\par
Un automorphisme $f$ est dit {\it de type parabolique} 
si le rayon spectral de $f^*$ sur 
$\H^{1,1}(V,\R)$
est \'egal \`a 1 et si
$f^*$ n'est pas d'ordre fini. Observons qu'alors $f$ est
d'entropie nulle et tout \'el\'ement de sa classe dans
$\mbox{Aut}^\#(V)$ est de type parabolique. On dira aussi que la classe de $f$ est 
{\it de type parabolique}.
L'argument donn\'e dans la
preuve du corollaire 2.2 montre que
sur $\H^{1,1}(V,\R)$ les valeurs propres de $f^*$ sont des racines de
l'unit\'e. En particulier,
il existe $n\geq 1$ tel que
sur $\H^{1,1}(V,\R)$, $(f^n)^*$ soit d\'efini par une matrice non
diagonalisable dont les valeurs propres sont \'egales \`a 1.  
\par
On dira qu'un groupe d'automorphismes $\G$ d'une vari\'et\'e compacte
complexe $V$ {\it pr\'eserve une classe} $c$ si pour tout $f\in \G$ 
il existe une constante $\lambda(f)$ telle que
$f^*(c)=\lambda(f)c$. On dira que la classe $c$ est {\it invariante} par
le groupe $\G$ et que $\lambda$ est {\it le coefficient multiplicatif}
associ\'e \`a $c$.
Notons $\langle f\rangle$ le groupe engendr\'e par $f$.
On a le lemme suivant:
\begin{lemme} Soit $\G$ un groupe commutatif 
d'automorphismes de $V$. Supposons que $\G$ contient un
  \'el\'ement $f$ d'entropie positive. Alors il existe 
une classe $c\in\overline \K$, $c\not=0$  
  qui est invariante par $\G$ et v\'erifie $f^* c=d_1(f)c$. 
\end{lemme}
\begin{preuve} Soit $\Gamma_1\subset \overline \K$ le c\^one des classes
  $c$ v\'erifiant $f^*c=d_1(f)c$. Puisque $\overline \K$ est un c\^one
convexe ferm\'e invariant par $f^*$, d'apr\`es une extension 
du th\'eor\`eme de 
Perron-Frobenius,
$\Gamma_1$ un c\^one convexe ferm\'e et non r\'eduit \`a
0. Il suffit pour cela d'utiliser la d\'ecomposition de Jordan de la matrice 
associ\'ee \`a $f^*$.
Soient $g\in \G$ et $c\in\Gamma_1$. On a 
$$f^*(g^*c)=g^*(f^*c)=d_1(f)g^*c.$$
Par cons\'equent, $g$ pr\'eserve $\Gamma_1$. 
Il existe $c\in \Gamma_1\setminus\{0\}$ telle que
$g^*c=\lambda(g)c$ o\`u $\lambda(g)$ est une constante positive. 
Notons $\Gamma_2\subset \Gamma_1$ le c\^one des classes $c$ v\'erifiant 
la derni\`ere relation. C'est un c\^one convexe ferm\'e.
Par induction en $g\in \G$, on construit une 
classe $c\in\Gamma_1$ invariante par 
$\G$. 
\end{preuve}
\par
On dira qu'un groupe $\G\subset \mbox{Aut}(V)$ est
{\it d'entropie positive} 
si tout \'el\'ement $f\in\G$ diff\'erent de l'identit\'e 
est d'entropie positive. Ce sont des groupes qui ne
contiennent pas d'\'el\'ements virtuellement isotopes \`a l'identit\'e
ni d'\'el\'ements de type parabolique, ils sont discrets.
\par
Pour tout $\tau=\{\tau(f)\}_{f\in \G}\in \R^\G$, notons
$\Gamma_\tau$ 
le c\^one des classes $c\in\overline\K$
telles que
$f^*(c)=\exp(\tau(f))c$ pour tout $f\in \G$. 
Si $\Gamma_\tau$ n'est pas r\'eduit \`a
$\{0\}$, $\exp(\tau(f))$ est une valeur propre de $f^*$ sur $\H^{1,1}(V,\R)$.
Notons $\F$ la famille des $\tau\in \R^\G$ tels que
$\Gamma_\tau\not=\{0\}$ et $\tau\not \equiv 0$. 
On montrera que $\F$ est fini et
on notera $\tau_1$, $\ldots$, $\tau_m$ les \'el\'ements de $\F$.
Notons \'egalement $\pi$ le morphisme de groupes
qui associe \`a $f\in \G$ le vecteur $(\tau_1(f),\ldots, \tau_m(f))\in \R^m$.
\begin{proposition} Le cardinal $m$ de $\F$ v\'erifie $m\leq h_1$ o\`u
  $h_1:=\dim \H^{1,1}(V,\R)$. 
Si $\G$ est commutatif 
d'entropie positive, 
$\pi$ est injective et 
son image est discr\`ete. 
De plus, $\G$ est commutatif libre et 
son rang $r$ v\'erifie $r\leq m$.
%
%
\end{proposition} 
\begin{preuve} Montrons que $m\leq h_1$, en particulier, $m$ est fini.
Puisque les $\tau_i$, pour tout $1\leq i<j\leq m$, sont
  diff\'erents, l'image $\pi(\G)$ de $\G$ n'est pas contenue dans
  l'hyperplan $(x_i=x_j)$. Par
  cons\'equent, le groupe additif $\pi(\G)$ n'est pas contenu dans la
  r\'eunion de ces hyperplans. Il existe donc un \'el\'ement $g_0\in \G$
  tel que $\tau_i(g_0)\not=\tau_j(g_0)$ pour tous $1\leq i<j\leq
  m$. Sur $\H^{1,1}(V,\R)$, l'application $g_0^*$ poss\`ede donc au moins 
$m$ valeurs propres distinctes. Par suite, $m\leq h_1$.
\par
Supposons que $\G$ soit commutatif et d'entropie positive.
D'apr\`es le lemme 4.1,
l'une des coordonn\'ees de $\pi(f)$ est \'egale \`a $\log d_1(f)$. Le
corollaire 2.2 implique que $\log d_1(f)$ est minor\'e
ind\'ependement de $f$. Par suite, $0$ est isol\'e dans 
$\pi(\G)$ et donc $\pi(\G)$ est discret. Comme $\G$ est
d'entropie positive, l'application $\pi$ est injective.
\par
Puisque $\pi(\G)$ est un sous groupe additif discret de $\R^m$, on a
$r\leq m\leq h_1$. 
\end{preuve}
\begin{lemme} Soient $c$, $c'$ et $c_i$ des classes dans $\overline\K$, 
$1\leq i\leq m\leq k-2$. Soit $g$ un automorphisme de $V$.
Supposons qu'il existe des constantes r\'eelles 
positives diff\'erentes $\lambda$ et $\lambda'$ telles que
$g^*(c_1\wedge\ldots\wedge c_m\wedge c)=\lambda c_1\wedge\ldots\wedge 
c_m\wedge c$ et 
$g^*(c_1\wedge\ldots\wedge c_m\wedge c')=
\lambda' c_1\wedge\ldots\wedge c_m\wedge c'$. Supposons de plus que 
$c_1\wedge\ldots\wedge c_m\wedge c\not=0$ et 
$c_1\wedge \ldots\wedge c_m\wedge c\wedge c'=0$.
Alors $c_1\wedge \ldots\wedge c_m\wedge c'=0$.
\end{lemme}
\begin{preuve}
D'apr\`es le corollaire
3.5, il existe $(a,b)\in \R^2\setminus\{0\}$ unique \`a une
constante multiplicative pr\`es, tel que 
$$c_1\wedge \ldots \wedge c_m \wedge (ac + b c')\wedge
c'_1 \wedge \ldots \wedge c'_{k-n-2} =0$$
pour toutes classes $c'_1$, $\ldots$, $c'_{k-m-2}$ dans $\H^{1,1}(V,\R)$. 
\par
Utilisant $g^*$, la dern\i\`ere \'egalit\'e implique que
$$c_1\wedge \ldots \wedge c_m \wedge (a\lambda c + b\lambda' c')\wedge
g^*c'_1 \wedge \ldots \wedge g^*c'_{k-m-2} =0.$$ 
Comme l'op\'erateur $g^*$ est inversible, on a 
$$c_1\wedge \ldots \wedge c_m \wedge (a\lambda c + 
b\lambda' c')\wedge
c'_1 \wedge \ldots \wedge c'_{k-m-2} =0$$
pour toutes classes $c'_1$, $\ldots$, $c'_{k-m-2}$ dans $\H^{1,1}(V,\R)$. 
Les couples $(a,b)$ et $(a\lambda,b\lambda')$ 
sont donc colin\'eaires. Comme $\lambda\not=\lambda'$, on a $a=0$
ou $b=0$. Or $c_1\wedge\ldots\wedge c_m\wedge c\not=0$, donc 
$a=0$. En prenant les classes $c'_i$ dans le
c\^one de K\"ahler, la derni\`ere \'egalit\'e implique que $c_1\wedge
\ldots \wedge c_m\wedge c'=0$ car les $c_i$ et $c'$  
appartiennent \`a $\overline \K$. 
\end{preuve}
\par
La proposition suivante implique le th\'eor\`eme principal. En effet, 
on a $h_k=1$ et la proposition (appliqu\'ee \`a $n=r=k$) implique que $r\not=k$. 
La propostion donne \'egalement des
  contraintes topologiques pour qu'une vari\'et\'e k\"ahl\'erienne
  poss\`ede un gros groupe commutatif 
d'automorphismes d'entropie positive.
\begin{proposition} Soit $\G$ un groupe commutatif d'automorphismes
  d'une vari\'et\'e k\"ahl\'erienne compacte $V$ de dimension $k\geq
  2$. Supposons que $\G$ soit d'entropie positive. 
Alors $\G$ est commutatif libre. Le rang $r$ de
  $\G$ v\'erifie ${r\choose n}\leq h_n$ pour tout $1\leq n\leq r$ 
o\`u $h_n:=\dim \H^{n,n}(V,\R)$; 
si $n$ divise $r$ on a ${r\choose n}\leq h_n-1$.
De plus, il existe $r+1$ classes
$c_1$, $\ldots$, $c_{r+1}$ dans $\overline \K$ telles que
  $c_1\wedge \ldots\wedge c_{r+1}\not=0$. 
\end{proposition}
\begin{preuve} L'application $\pi$ \'etant injective, $\pi(\G)$
  engendre un sous-espace de dimension $r$. On peut supposer que la projection de ce
  sous-espace sur les $r$ premiers coordonn\'ees est bijective. Cela
  signifie que $\tau_1$, $\ldots$, $\tau_r$ est une famille libre engendrant
  $\F$. 
Notons $\Pi$ le morphisme qui associe \`a $g\in \G$ le vecteur
  $(\tau_1(g),\ldots,\tau_r(g))\in \R^r$. L'application $\Pi$ est donc
  injective et son image engendre l'espace $\R^r$.
Soit $c_i$ une classe non nulle dans $\Gamma_{\tau_i}$, $1\leq
  i\leq r$.
\par
Pour tout $I=(i_1,\ldots, i_n)\subset \{1,\ldots,r\}$, notons
$c_I:=c_{i_1}\wedge \ldots \wedge c_{i_n}$ et $\tau_I:=
\tau_{i_1}+\cdots +\tau_{i_n}$. Montrons que $c_I\not =0$ quand
$|I|\leq k$. Sinon, on peut supposer que
$c_1\wedge \ldots \wedge c_{n-1}\not =0$, 
$c_1\wedge \ldots \wedge c_{n-2}\wedge c_n\not =0$
 et $c_1\wedge \ldots
\wedge c_n=0$ avec $2\leq n\leq \min(k,r)$. On verra que $\min(k,r)$
est toujours \'egal \`a $r$. 
Soit $g\in \G$. On peut appliquer 
lemme 4.3 pour $m:=n-2$, $c:=c_{n-1}$ et $c':=c_n$. On constate que 
$\tau_{n-1}(g)= \tau_n(g)$. Ceci implique que $\Pi(\G)$ est contenue
dans l'hyperplan $\{x_{n-1}=x_n\}$. C'est contraire au choix de $\Pi$.
\par
Montrons que les classes $c_I$ avec $|I|=n$ 
sont lin\'eairement ind\'ependantes. Sinon, il existe des nombres
r\'eels $a_I$ non tous nuls, tels que
$$\sum a_I c_I=0.$$
Le fait que les $c_i$  soient invariantes par $g^p$ pour tout $g\in \G$ et
tout $p\in \Z$,  implique que 
$$\sum a_I \exp(p\tau_I(g))c_I=0.$$
Par cons\'equent, il existe $I\not=J$ tels que
$\tau_I(g)=\tau_J(g)$. Le groupe $\Pi(\G)$ est donc contenu dans la
r\'eunion finie des hyperplans du type
$$\sum_{i\in I} x_i = \sum_{i\in J} x_j.$$
Par suite, il est contenu dans un hyperplan. 
C'est contraire au choix de $r$. 
\par
Les $r\choose n$ classes $c_I$ sont lin\'eairement ind\'ependantes
dans $\H^{n,n}(V,\R)$. On en d\'eduit que ${r\choose n}\leq h_n$ pour
tout $n\leq \min(k,r)$. Si
$r\geq k$, en prenant $n=k$, la derni\`ere in\'egalit\'e donne $r\leq
k$ et donc $\min(k,r)=r$. 
\par
Supposons maintenant que $n$ divise $r$. Posons $r=mn$. Montrons que 
 ${r\choose n}\not = h_n$. Sinon, les classes $c_I$ forment une base de 
$\H^{n,n}(V,\R)$ car elles sont lin\'eairement ind\'ependantes. 
On peut alors choisir $m$ classes $c_{I_1}$, $\ldots$, $c_{I_m}$ telles que
$c_{I_1}\wedge\ldots\wedge c_{I_m}\not =0$. La derni\`ere classe appartient 
\`a $\H^{k,k}(V,\R)$. Par cons\'equent, 
$g^* (c_{I_1}\wedge\ldots\wedge c_{I_m})= c_{I_1}\wedge\ldots\wedge c_{I_m}$
pour tout $g\in \G$. On en d\'eduit que $\sum \tau_{I_i}=0$.
Par cons\'equent, $\Pi(\G)$ est contenue dans un hyperplan, 
ce qui n'est pas possible. On a montr\'e que ${r\choose n}\leq h_n-1$. 
Ceci (pour $r=k=n$) implique que $\Aut(V)$ ne contient pas de sous-groupe
commutatif, libre, de rang $k$ et d'entropie positive.
\par
Il reste \`a construire une classe $c_{r+1}\in \overline \K$ telle que 
$c_1\wedge\ldots\wedge c_{r+1}\not =0$. 
On sait d\'ej\`a que $c_1\wedge\ldots\wedge c_r\not =0$.
Puisque $\Pi(\G)$ est un r\'eseau 
qui engendre $\R^r$, il existe $f\in \G$ tel que les coordonn\'ees de 
$\Pi(f)$ soient toutes strictement n\'egatives. D'apr\`es le lemme 4.1,
il existe une classe $c_{r+1}\in\overline \K$ invariante par $\G$ 
telle que $f^*c_{r+1}=d_1(f) c_{r+1}$. D'apr\`es le lemme 4.3, on a 
$c_1\wedge\ldots\wedge c_{r+1}\not =0$ car sinon, on en d\'eduit que 
$c_1\wedge\ldots\wedge c_r=0$.
\end{preuve}
\begin{exemple} \rm
Soit $\mathbb{T}:=\C/(\Z+i\Z)$ le tore complexe associ\'e au r\'eseau
$\Z+i\Z$. Posons $V:=\mathbb{T}^k$. Notons $\mbox{SL}(k,\Z)$ le groupe
des matrices de rang $k$, de d\'eterminant $1$ \`a coefficients
entiers. 
Ce groupe agit sur $V$ et contient
des sous-groupes commutatifs libres de rang $k-1$, diagonalisables
dans $\R$ donc d'entropie positive. 
En effet, le groupe $\Delta$ des matrices r\'eelles diagonales 
agit \`a gauche sur l'espace homog\`ene
$\mbox{SL}(k,\R)/\mbox{SL}(k,\Z)$ et poss\`ede des orbites
compactes (\voir par exemple Samuel \cite[p.72]{Samuel}
ou Prasad-Raghunathan \cite[theorem
2.8]{PrasadRaghunathan}). 
Si l'orbite de $g$ est compacte, $g^{-1}\Delta g\cap
\mbox{SL}(k,\Z)$ est un groupe commutatif libre de rang $k-1$.
Cet exemple montre que la borne sup\'erieure $k-1$ trouv\'ee dans le
th\'eor\`eme 4.6 est optimale. 
\end{exemple}
\begin{theoreme} 
Soit $\G$ un groupe d'automorphismes
  d'une vari\'et\'e k\"ahl\'erienne compacte $V$ de dimension $k\geq
  2$. Supposons que $\G$ soit d'entropie positive. 
Supposons \'egalement que $\G$ pr\'eserve $k-1$ classes $c_1$, $\ldots$,
  $c_{k-1}$ de $\overline \K$ qui v\'erifient $c_1\wedge \ldots\wedge
  c_{k-1}\not=0$. Alors $\G$ est commutatif, libre et de rang $r\leq k-1$. 
\end{theoreme}
\begin{preuve}
Comme pr\'ec\'edemment, on note $\tau_i$ la
fonction associ\'ee \`a $c_i$ et $\pi:\G\longrightarrow \R^{k-1}$
l'application associ\'ee aux classes invariantes $c_1$, $\ldots$, $c_{k-1}$.
Montrons que $\pi$ est injective.
Supposons qu'il existe $g\in \G$ tel que $g^*c_i=c_i$ pour 
tout $1\leq i\leq k-1$. 
Soit $c$ une classe de $\overline \K\setminus\{0\}$ 
v\'erifiant $g^*c=\lambda c$ o\`u $\lambda:=d_1(g)$. On a 
$g^*(c\wedge c_1\wedge 
\ldots \wedge c_{k-1}) =\lambda c\wedge c_1\wedge 
\ldots \wedge c_{k-1}$. Or $g^*$ est \'egal \`a l'identit\'e 
sur $\H^{k,k}(V,\R)$ et 
$\lambda>1$. Donc $c\wedge c_1\wedge \ldots\wedge c_{k-1}=0$. 
Utilisant le lemme 4.3, on montre par r\'ecurrence que $c=0$, ce qui est 
impossible. Donc $\pi$ est injective et $\G$ est commutatif.
D'apr\`es le th\'eor\`eme principal, on a $r\leq k-1$.
\end{preuve}
\begin{theoreme} Soit $\G'$ un groupe commutatif d'automorphismes d'une
vari\'et\'e k\"ahl\'erienne compacte $V$ de dimension $k\geq 2$. 
Alors l'ensemble  $U$ des \'el\'ements d'entropie nulle 
de $\G'$ est un groupe. De plus, 
il existe un sous-groupe commutatif libre, d'entropie 
positive $\G$ de $\G'$, de 
rang $r\leq k-1$ tel que $\G'\simeq U\times \G$. 
Si $r=k-1$, le groupe $U$ 
est fini.
\end{theoreme}
\begin{preuve} Comme pr\'ec\'edemment, on construit $r+1$ classes $c_1$, 
$\ldots$, $c_{r+1}$ de $\overline \K$, invariantes par $\G'$ et v\'erifiant
$c_1\wedge\ldots\wedge c_{r+1}\not=0$. 
Les $r+1$ fonctions $\tau_i$ associ\'ees 
\`a ces classes sont diff\'erentes. On suppose ici que $r$ est le nombre 
maximal tel qu'il existe $r$ fonctions $\tau_i$ 
lin\'eairement ind\'ependantes. 
Notons toujours $\pi:\G'\longrightarrow \R^{r+1}$ l'application 
associ\'ee aux $\tau_i$. L'image de $\pi$ est un r\'eseau qui engendre
un hyperplan de $\R^{r+1}$. Seuls les \'el\'ements d'entropie nulle
s'envoient \`a 0. On a $U=\pi^{-1}(0)$ et donc $U$ est un groupe.
On choisit des \'el\'ements $g_1$, $\ldots$, $g_r$ de $\G'$ tels 
que les $\pi(g_i)$  
engendrent le r\'eseau $\pi(\G')$. Notons $\G$ le groupe engendr\'e par 
$g_1$, $\ldots$, $g_r$. Alors $\G$ est commutatif libre de 
rang $r$ et d'entropie positive
car il est isomorphe 
\`a $\pi(\G')=\pi(\G)$. Il est clair que $\G'\simeq U\times \G$. 
\par
Supposons maintenant que $r=k-1$.
Posons $c:=c_1+\cdots+c_k$. 
Pour tout $u\in U$, on a $u^*c_i=c_i$ et donc $u^*c=c$.
Puisque $c_1\wedge\ldots\wedge c_k\not=0$, on a $c^k\not =0$ et $\pi(\G')$ 
est contenu dans l'hyperplan $\{x_1+\cdots+x_k=0\}$. 
D'apr\`es un th\'eor\`eme de
Demailly-Paun \cite{DemaillyPaun}, le c\^one de K\"ahler 
est une des composantes connexes
de l'ouvert des classes $[\alpha]$ telles que $\int_V \alpha^k>0$. 
Or $c\in \overline \K$ et $c^k\not=0$, donc $c\in\K$.
On a $U\subset \Aut_c(V)$. 
Comme on l'a rappel\'e, d'apr\`es Lieberman
\cite[Proposition 2.2]{Lieberman}, 
$\Aut_c(V)/\Aut(V)$ est fini. Posons $U_0:=U\cap\Aut_0(V)$. 
Il suffit de montrer que $U_0$ est fini.
\par
Rappelons que $\Aut_0(V)$ est compactifiable
\cite[Theorem 3.3]{Lieberman} et 
notons $\overline U_0$ l'adh\'erence de Zariski de $U_0$ dans
$\Aut_0(V)$. D'apr\`es Lieberman \cite[Corollary 3.6]{Lieberman},
$\overline U_0$ est commutatif. Posons $m:=\dim_\C \overline U_0$.
\par
Si $m\geq k$, $\overline U_0$ contient
un groupe de Lie commutatif
de dimension $k$ agissant sur $V$. On en d\'eduit que $V$ est un
tore complexe. Dans ce cas, si $z$ d\'esigne un point du tore,
tout automorphisme $g$ de $V$
s'\'ecrit $g(z)=A_gz+B_g$ o\`u $A_g$ est une matrice carr\'ee de rang
$k$, $B_g$ est un vecteur de longueur $k$. 
Pour un \'el\'ement $g$ de $\Aut_0(V)$,
on a $A_g=\id$.
Maintenant, fixons un $g\in \G$ tel que
$\tau_i(g)>0$ pour $i=1,\ldots, k-1$. Il est clair que toute
valeur propre de $A_g$ est de module diff\'erent de $1$. On
v\'erifie facilement qu'il n'y a qu'un nombre fini d'\'el\'ements de
$\Aut_0(V)$ qui commutent avec $g$. Ceci implique le th\'eor\`eme. 
\par
Supposons que $0<m\leq k-1$.
Soit $g\in\G$. D'apr\`es Lieberman \cite[Proposition 3.7]{Lieberman},
le groupe $g\overline U_0g^{-1}$ est Zariski ferm\'e. Comme ce
groupe contient $U_0$, il contient aussi $\overline U_0$. On en
d\'eduit que $g\overline U_0 g^{-1}=\overline U_0$ pour tout
$g\in\G$. Par
cons\'equent, $\G$ pr\'eserve la fibration engendr\'ee par les 
orbites compactifi\'ees de $\overline U_0$.
Soit
$\Omega$ la
classe associ\'ee \`a une orbite g\'en\'erique compactifi\'ee de
$\overline U_0$. Cette classe ne d\'epend pas de l'orbite 
g\'en\'erique choisie.
Comme $c=c_1+\cdots+c_k$ est une classe de K\"ahler, il existe 
$1\leq i_1\leq \cdots\leq i_m\leq k$ tels que $\Omega\wedge
c_{i_1} \wedge \ldots\wedge c_{i_m}\not =0$. Or $g^*$ est
l'identit\'e sur $\H^{k,k}(V,\R)$. On a
donc $\tau_{i_1}+\cdots+\tau_{i_m}=0$. Par cons\'equent, $\pi(\G)$ est
contenu dans l'hyperplan  $x_{i_1}+\cdots+x_{i_m}=0$ de $\R^k$. Ceci
contredit le fait que $\pi(\G)$ engendre l'hyperplan
$\{x_1+\cdots+x_k=0\}$.
\end{preuve}
\begin{theoreme} Soit $V$ une vari\'et\'e k\"ahl\'erienne compacte de
  dimension $k\geq 2$. Soit $\G^\#$ un sous-groupe  
de $\mbox{\rm Aut}^\#(V)$. Alors
\begin{enumerate}
\item Si $\G^\#$ est commutatif et d'entropie positive (\cad que $\G^\#$ 
ne contient pas de classes de type parabolique), alors $\G^\#$ est 
commutatif libre et son rang $r$ v\'erifie $r\leq k-1$.
\item Si $\G^\#$ pr\'eserve $k-1$ classes $c_1$, $\ldots$, $c_{k-1}$ dans
$\overline \K$ qui v\'erifient $c_1\wedge \ldots\wedge c_{k-1}\not =0$,
alors $\G^\#$ est commutatif, libre et d'entropie positive.
\item Si $\G^\#$ est commutatif et s'il contient un 
sous-groupe $\H^\#$ commutatif libre 
de rang $k-1$ et d'entropie positive, alors 
$\G^\#/\H^\#$ est fini.
\end{enumerate}
\end{theoreme}
\begin{preuve} Puisque $\mbox{Aut}_0(V)$ agit trivialement sur
  $\H^{1,1}(V,\R)$, on peut adapter les preuves des propositions 4.4, 4.7 et 
du th\'eor\`eme 4.6.
\end{preuve}
\begin{remarque} \rm
Si $\G^\#$ est un groupe commutatif contenant 
des classes de type parabolique, le rang de $\G^\#$
peut d\'epasser $k-1$. Soit en effet
$\mathbb{T}^2$ un tore de dimension 2
o\`u $\mathbb{T}:=\C/(\Z+ i\Z)$. Consid\'erons le groupe 
des automorphismes d\'efinis par
les matrices 
$$A_{m,n}:=
\left(\begin{array}{cc}
1 & m+in \\
0 & 1
\end{array}\right)$$
avec $m$, $n$ entiers.
C'est un groupe commutatif, libre et de rang $2$. Sur $\mathbb{T}^k$ on
peut contruire un groupe analogue de rang $2k-2$.
\end{remarque}
\begin{remarque}\rm
On dira que la vari\'et\'e $V$ est {\it multisym\'etrique} 
si elle admet un groupe 
d'automorphismes  
commutatif, libre, de rang $k-1$ et d'entropie positive.  
Il serait int\'eressant de classer les vari\'et\'es multisym\'etriques. 
En dimension 2, il s'agit de
classer les surfaces admettant au moins un automorphisme d'entropie
maximale. Le lecteur peut consulter les articles de Cantat 
\cite{Cantat1}, Mazur \cite{Mazur} et McMullen \cite{McMullen} pour une liste
d'exemples et une \'etude dynamique des automorphismes d'une surface
complexe. 
Les exemples de Mazur \cite{Mazur} 
peuvent \^etre \'etendus en dimension strictement sup\'erieure \`a 2.
\end{remarque}
\small
\par\noindent
Tien-Cuong Dinh et Nessim Sibony\\
Math\'ematique - B\^at. 425, UMR 8628\\
Universit\'e Paris-Sud, 91405 Orsay, France.\\ 
E-mails: Tiencuong.Dinh@math.u-psud.fr, Nessim.Sibony@math.u-psud.fr
\end{document}